\magnification1200
\overfullrule=0pt
\vsize=8truein
\hsize=6truein
\hoffset=18truept
\voffset=12truept

\font\sixbsy=cmbsy6 
\font\sixbf=cmbx6
\font\sixbfit=cmmib6 
\font\eightbsy=cmbsy8
\font\eightbf=cmbx8
\font\eightbfit=cmmib8
\font\twelvebf=cmbx10 at 12pt
\font\twelvebfit=cmmib10 at 12pt
\font\twelvebsy=cmbsy10 at 12pt

\font\tenmsb=msbm10 
\font\sevenmsb=msbm7
\font\fivemsb=msbm5
\newfam\msbfam
\textfont\msbfam=\tenmsb
\scriptfont\msbfam=\sevenmsb
\scriptscriptfont\msbfam=\fivemsb
\def\Bbb{\fam\msbfam\tenmsb}

\def\bn{\medskip \noindent}
\def\sn{\smallskip \noindent}
\def\O{{\cal O}}
\def\La{\longrightarrow}
\def\F{{\cal F}}
\def\H{{\cal H}}
\def\L{{\cal L}}
\def\Pic{{\rm Pic}}
\def\f{{\hat f}}
\def\X{{\hat X}}

\def\G{{\cal G}}
\def\C{{\cal C}}
\def\bC{{\Bbb C}}
\def\bN{{\Bbb N}}
\def\bP{{\Bbb P}}
\def\bQ{{\Bbb Q}}
\def\bR{{\Bbb R}}
\def\bZ{{\Bbb Z}}
\def\Ker{\mathop{\rm Ker}\nolimits}
\def\Im{\mathop{\rm Im}\nolimits}
\def\Sl{\mathop{\rm Sl}\nolimits}
\def\SU{\mathop{\rm SU}\nolimits}
\def\rk{\mathop{\rm rk}\nolimits}
\def\id{\mathop{\rm id}\nolimits}
\def\pr{\mathop{\rm pr}\nolimits}
\def\exp{\mathop{\rm exp}\nolimits}
\def\log{\mathop{\rm log}\nolimits}
\def\section#1{\noindent{%
\textfont0=\twelvebf \scriptfont0=\eightbf \scriptscriptfont0=\sixbf
\textfont1=\twelvebfit \scriptfont1=\eightbfit \scriptscriptfont1=\sixbfit
\textfont2=\twelvebsy \scriptfont2=\eightbsy \scriptscriptfont2=\sixbsy
\twelvebf #1}}
\def\dasharrow{\mathrel{\hbox{
\kern-1pt \vrule height2.45pt depth-2.15pt width2pt} \kern1pt{\vrule
height2.45pt depth-2.15pt width2pt} \kern1pt{\vrule height2.45pt
depth-2.15pt width2pt} \kern1pt{\vrule height2.45pt depth-2.15pt
width1.7pt\kern-1.7pt} {\raise1.4pt\hbox{$\scriptscriptstyle\succ$}}}}

\catcode`\@=11
\newskip\ttglue
\font\ninerm=cmr9

\font\sixrm=cmr6

\font\ninei=cmmi9
\font\eighti=cmmi8
\font\sixi=cmmi6
\skewchar\ninei='177 \skewchar\eighti='177 \skewchar\sixi='177

\font\ninesy=cmsy9
\font\eightsy=cmsy8
\font\sixsy=cmsy6
\skewchar\ninesy='60 \skewchar\eightsy='60 \skewchar\sixsy='60

\font\ninebf=cmbx9
\font\eightbf=cmbx8
\font\sixbf=cmbx6

\font\ninett=cmtt9
\font\eighttt=cmtt8

\hyphenchar\tentt=-1 
\hyphenchar\ninett=-1
\hyphenchar\eighttt=-1

\font\ninesl=cmsl9

\font\nineit=cmti9

\def\ninepoint{\def\rm{\fam0\ninerm}%
  \textfont0=\ninerm \scriptfont0=\sixrm \scriptscriptfont0=\fiverm
  \textfont1=\ninei \scriptfont1=\sixi \scriptscriptfont1=\fivei
  \textfont2=\ninesy \scriptfont2=\sixsy \scriptscriptfont2=\fivesy
  \textfont3=\tenex \scriptfont3=\tenex \scriptscriptfont3=\tenex
  \def\it{\fam\itfam\nineit}%
  \textfont\itfam=\nineit
  \def\sl{\fam\slfam\ninesl}%
  \textfont\slfam=\ninesl
  \def\bf{\fam\bffam\ninebf}%
  \textfont\bffam=\ninebf \scriptfont\bffam=\sixbf
   \scriptscriptfont\bffam=\fivebf
  \def\tt{\fam\ttfam\ninett}%
  \textfont\ttfam=\ninett
  \tt \ttglue=.5em plus.25em minus.15em
  \normalbaselineskip=11pt
  \def\MF{{\manual hijk}\-{\manual lmnj}}%
  \let\sc=\sevenrm
  \let\big=\ninebig
  \setbox\strutbox=\hbox{\vrule height8pt depth3pt width\z@}%
  \normalbaselines\rm}

\font\smaller=cmr8
\font\ninecsc=cmcsc9
\font\tencsc=cmcsc10

\def\section#1{\bigskip\medskip
\centerline{\tencsc #1}\medskip}

\headline={\ifnum\pageno>1\smaller\ifodd\pageno\hfill 
THE ALGEBRAIC DIMENSION OF COMPACT COMPLEX THREEFOLDS
\hfill\the\pageno \else
\the\pageno\hfill\uppercase{
Fr\'ed\'eric Campana, Jean-Pierre Demailly, and
Thomas Peternell}\hfill\fi\else\hss\fi}

\footline={\hss}

\centerline{\bf \uppercase{The Algebraic Dimension of Compact Complex Threefolds}}
\smallskip
\centerline{\bf \uppercase{with vanishing second Betti Number}}
\footnote{}{\ninepoint 
\noindent Parts of this work were done during one of the author's
stay at MSRI. Research at MSRI is supported in part by NSF grant DMS-9022140.}
\bigskip\medskip
\centerline{\ninepoint\uppercase{
Fr\'ed\'eric Campana, Jean-Pierre Demailly, and
Thomas Peternell}}

\midinsert
\narrower\narrower
\noindent {\ninecsc Abstract.}
\baselineskip=10pt
{\ninepoint 
We investigate compact complex manifolds of
dimension three and second Betti number $b_2(X) = 0$.
We are interested in the
algebraic dimension $a(X)$, which is by definition the transcendence degree
of the field of meromorphic functions over the field of complex numbers.
The topological Euler characteristic $\chi_{\rm
top}(X) $ equals the third Chern class $c_3(X)$ by a theorem of
Hopf. Our main result is that, if
$X$ is a compact 3-dimensional complex manifold with $b_2(X) =
0$ and $a(X) > 0$, then $c_3(X) = \chi_{\rm top}(X) = 0$, that is, we either
have $b_1(X) = 0, \ b_3(X) = 2$ or $b_1(X) = 1, \ b_3(X) = 0.$ 

}
\endinsert

\section{Introduction}

\bn In this note we shall investigate compact complex manifolds of
dimension three and second Betti number $b_2(X) = 0.$ Of course such a
manifold cannot be algebraic. Therefore we will be interested in the
algebraic dimension $a(X)$ which is by definition the transcendence degree
of the field of meromorphic functions over the field of complex numbers.
Note that $a(X) > 0$ if and only if $X$ admits a non-constant meromorphic
function. The topological Euler characteristic will be denoted $\chi_{\rm
top}(X) $ which is also the third Chern class $c_3(X)$ by a theorem of
Hopf. Our main result is

\bn {\bf Theorem}
\sn {\it Let $X$ be a compact 3-dimensional complex manifold with $b_2(X) =
0$ and $a(X) > 0.$ Then $c_3(X) = \chi_{\rm top}(X) = 0,$ i.e. we either
have $b_1(X) = 0, \ b_3(X) = 2$ or $b_1(X) = 1, \ b_3(X) = 0.$ }

\bn Notice that if $a(X) = 3,$ i.e. $X$ is Moishezon, then we have $b_2(X)
> 0,$ but examples of compact threefolds $X$ with $a(X) = 1$ or 2 and with
the above Betti numbers exist (see sect. 4). We also show by an example
that the assumption $a(X) > 0$ cannot be omitted.

\bn The following corollary was actually our motivation for the Theorem:

\bn {\bf Corollary}\par\nobreak
\sn {\it Let $X$ be a compact complex manifold homeomorphic to the
6-dimensional sphere $S^6.$ Then $a(X) = 0.$ }

\bn In other words, $S^6$ does not admit a complex structure with a
non-constant meromorphic function.

\bn The proof of the Theorem consists in showing that for generic $\L \in
{\Pic}^0(X),$ we have
$$
h^0(X,B \otimes \L) = h^2(X,B \otimes \L)= 0 $$
for $B$ either the tangent bundle $T_X$ or the cotangent bundle
$\Omega^1_X.$ This forces ${c_3(X) \over 2} = \chi(X,T_X) $ to be 0. If
$a(X) = 2$ or if $a(X) = 1$ and the algebraic reduction is not holomorphic,
we can even prove the last vanishing statements for all vector bundles $B$
on~$X.$

\sn In the last section we study more closely the structure of threefolds
$X$ with $b_2(X) = 0$ and algebraic dimension 1 whose algebraic reduction
is holomorphic. We show e.g. that smooth fibers can only be Inoue surfaces,
Hopf surface with algebraic dimension 0 or tori.

\section{1. Preliminaries and Criteria for the vanishing of $H^0.$} \bn
{\bf 1.0 Notations}

\sn (1) Let $X$ be a compact complex manifold, always assumed to be
connected. The algebraic dimension, denoted $a(X),$ is the transcendence
degree of the field of meromorphic functions over $\bC.$ \sn (2) $b_i(X) =
\dim H^{i}(X,\bR)$ denotes the $i-$th Betti number of $X.$ \sn (3) If $G$
is a finitely generated abelian group, then $rk G$ will denote its rank
(over $\bZ).$ \sn (4) If $X$ is a compact space, then $h^q(X,\F)$ denotes
the dimension of $H^q(X,\F).$

\bn {\bf 1.1 Proposition}
\sn {\it Let $Y$ be a connected compact complex space (not necessarily
reduced), every component $Y_i$ of $Y$ being of positive dimension. Let $D$
be an effective divisor on $Y$ such that $D \vert Y_i \ne 0$ for all $i.$
Let $\F$ be a locally free sheaf on $Y.$ Then there exists $k_0 \in \bN$
such that
$$
H^0(Y,\F \otimes \O_Y(-kD)) = 0
$$
for $k \geq k_0.$}

\bn {\bf Proof.} We have natural inclusions $$
H^0(Y,\F \otimes \O_Y(-(k+1)D)) \subset H^0(Y;\F \otimes \O_Y(-kD)). $$
Take $k_0$ such that this sequence is stationary for $k \geq k_0.$ Then $s$
has to vanish at any order along $D\vert Y_i$ for every $i$ ($s$ can be
thought of locally as a tuple of holomorphic functions), hence $s \vert Y_i
= 0,$ and $s = 0.$

\bn {\bf 1.2 Corollary}
\sn {\it Let $S$ be a smooth compact complex surface containing an
effective divisor $C$ such that $c_1(\O_S(C)) = 0.$ Let $B$ be a vector
bundle on~$S.$ Then for a generic $\L \in \Pic^0(X)$ we have $$
H^0(S,B \otimes \L) = H^2(S,B \otimes \L) = 0. $$
In particular $\chi(S,B) = -h^1(S,B \otimes \L) \leq 0.$}

\bn {\bf Proof.} The vanishing $H^0(S,B \otimes \O_S(-kC)) = 0$ for large
$k$ follows from (1.1). Since $\O_S(kC)$ is topologically trivial for large
suitable $k,$ the required $H^0-$vanishing follows from semi-continuity.
The $H^2-$vanishing follows by applying the previous arguments to $B^*
\otimes K_S$ and Serre duality.

\bn {\bf 1.3 Corollary}
\sn {\it Let $X$ be a smooth compact threefold with $b_2(X) = 0$ carrying
an effective divisor $D.$ Then $H^0(X,B \otimes \L) = 0$ for generic $\L$
in $\Pic^0(X)$ and every vector bundle $B$ on~$X.$ }

\bn {\bf Proof.} Since $c_1(\O_X(mD)) = 0$ in $H^2(X,\bZ)$, we can apply
(1.1) to obtain
$$
H^0(X,B \otimes \O_X(-kmD)) = 0
$$
for $k$ large. Now we conclude again by semi-continuity.

\bn The last lemma is of course well-known; we include it for the
convenience of the reader.

\bn {\bf 1.4 Lemma}
\sn {\it Let $X$ be a compact manifold of dimension $n$ with $a(X) = n.$
Then $b_2(X) > 0.$}

\bn {\bf Proof.} Choose a birational morphism $\pi : \hat X \La X$ such
that $\hat X$ is a projective manifold. Take a general very ample divisor
$\hat D$ on $\hat X$ and a general curve $\hat C \in \hat X.$ Let $ D = \pi
(\hat D)$ and $C = \pi(\hat C).$ Then $D$ meets $C$ in finitely many
points, hence $D \cdot C > 0,$ in particular $c_1(\O_X(D)) $ is not torsion
in $H^2(X;\bZ).$

\bn {\bf 1.5 Lemma}\par\nobreak
\sn {\it Let $X$ be a smooth compact threefold and $f : X \La C$ be a
surjective holomorphic map to the smooth curve $C$. Let $\F$ be a locally
free sheaf on~$X.$ Then $R^{i}f_*(\F)$ is locally free for all $i.$}

\bn {\bf Proof.} (a) Note that local freeness is equivalent to torsion
freeness, since dim$C = 1.$ Hence the claim is clear for $i = 0.$

\sn (b) Next we treat the case $ i = 2.$ We shall use relative duality (see
[RRV71], [We85]); it states in our special situation ($f $ is flat with
even Gorenstein fibers) thatif $R^jf_*(\G)$ is locally free for a given
locally free sheaf $\G$ and fixed $j,$ then $$ R^{2-j}f_*(\G^* \otimes
\omega_{X\vert C}) \simeq R^jf_*(\G)^*,$$ in particular $R^{2-j}f_*(\G^*
\otimes \omega_{X \vert C}) $ is locally free. Here $\omega_{X \vert C} =
\omega_X \otimes f^*(\omega_C^*) $ is the relative dualising sheaf.
Applying this to $j = 0$ and $\G = \F^* \otimes \omega^*_{X \vert C}$ our
claim for $i = 2$ follows.

\sn (c) Finally we prove the freeness of $R^1f_*(\F).$ By a standard
theorem of Grauert it is sufficient that $h^1(X_y,\F \vert X_y)$ is
constant, $X_y$ the analytic fiber over $y \in C.$ By flatness,
$\chi(X_y,\F \vert X_y)$ is constant, hence it is sufficient that
$h^{j}(X_y,\F \vert X_y)$ is constant for $j=0$ and $j=2.$ By the vanishing
$R^3f_*(\F) = 0$, we have (see e.g. [BaSt76]) $$
R^2f_*(\F)\vert \{ y\} \simeq H^2(X_y, \F \vert X_y). $$
Therefore $h^2(X_y,\F \vert X_y)$ is constant by (b). Finally $$
h^{0}(X_y, \F \vert X_y) = h^2(X_y, \F^* \vert X_y \otimes \omega_{X_y}) =
h^2(X_y, \F^* \otimes \omega_X \vert X_y) $$
is constant by applying the same argument, now to $\F^* \otimes \omega_X.$

\section{2. Criteria for the vanishing of $H^2(X,B).$}

\bn {\bf 2.0 Set-up}

\sn We fix in this section the following situation: \sn (1) $X$ is a smooth
compact threefold with $b_2(X) = 0.$ \sn (2) $V$ is a normal projective
variety of dimension $d = 1$ or $2$ \sn (3) $f : X \dasharrow V$ is a
meromorphic map \sn (4) $\hat X$ is a normal compact threefold with a
birational map $\sigma : \hat X \La X$ such that the induced map $\hat f :
\hat X \La V$ is holomorphic; $\hat X$ is smooth in case $d = 1$ \sn (5)
$\hat f$ is equidimensional.

\bn In our applications $f$ will be an algebraic reduction.

\bn {\bf 2.1 Lemma}
\sn {\it Let $f : X \La V$ be as in 2.0,(1)-(3). Then (4) and (5) can
always be achieved, possibly after changing $V$ birationally.}

\bn {\bf Proof.} The claim is obvious if $d=1$, by elimination of
indeterminacies. So let $d = 2.$ Take a sequence of blow-ups $\rho : X_1
\La X $ such that the induced map $f_1 : X_1 \La V$ is holomorphic. Then
however $f_1$ will in general not be equidimensional. To achieve this we
could either apply Hironaka's flattening theorem or we can argue in a more
elementary way as follows. Let $V_0 = \{y \in V \vert f_1^{-1}(y) {\rm \ is
\ 1-dimensional} \}$. We obtain a holomorphic map $\mu : V_0 \La
\C_1(X_1),$ the space of 1-cycles on $X_1.$ It is a standard fact in the
theory of cycle spaces that $\mu$ extends to a meromorphic map $\mu : V
\dasharrow \C_1(X_1).$ Now let $V' = {\overline {\mu(V)}} \subset
\C_1(X_1)$ and let $X' \La V'$ be the induced family. Modify $V'$ to $V''$
so that $V''$ is projective; then take $X'' = X' \times _{V'} V''$ (base
change). Finally let $\hat X$ be the normalisation of $X''$ and $\hat V$
the normalisation of $V''.$ Then the induced map $\hat f : \hat X \La \hat
V$ has the properties we are looking for. \bn From now on we substitute
$\hat V$ by $V.$

\bn First we will treat the case $d = 2.$ We fix an ample effective divisor
$A$ on $V$ in the set-up (2.0) and let $\L = \sigma_*(\f^*(A))^{**}.$ Then
$$
\f^*(A) = \sigma^*(\L) \otimes \O_{\X}(-E) $$
with an effective divisor $E$ on $\X$ supported on the exceptional locus of
$\sigma.$

\bn {\bf 2.2 Proposition}
\sn {\it $H^2(X,B \otimes \L^{\otimes m}) = 0$ for every vector bundle $B$
and $m \gg 0$ (depending on $B.$)}

\bn {\bf Proof.} First we claim
$$
H^2(\X,\sigma^*(B \otimes \L^{\otimes m})(-mE)) = 0 \eqno (1) $$
for large $m.$

\sn For the proof notice $\sigma^*(B \otimes \L^{\otimes m}) (-mE)) \simeq
\sigma^*(B) \otimes \f^*(mA).$ We consider the Leray spectral sequence
associated to $\f.$ Then
$$
E_2^{2,0} = H^2(V,\f_*(\sigma^*(B)) \otimes mA) = 0, m >> 0,$$ and $$
E^{1,1}_2 = H^1(V,R^1\f_*(\sigma^*(B)) \otimes mA) = 0, m >> 0 $$
by Serre's vanishing theorem. Moreover $E_2^{0,2} = 0$ simply because
$R^2\f_*(\sigma^*(B)) = 0,$ all fibers of $\f$ being 1-dimensional. Now the
claim (1) follows from the spectral sequence.

\sn Let $A_m$ be the complex subspace of $X$ defined by the ideal sheaf
$\sigma_*(-mE).$ Then we consider the following commutative diagram, the
vertical arrows just being pull-back maps (as arising in the Leray spectral
sequence)
{\ninepoint
$$
\matrix{
H^2(\X,\sigma^*(B \otimes \L^m)(-mE)) & \La & H^2(\X, \sigma^*(B \otimes
\L^m)) & \La & H^2(mE,\sigma^*(B \otimes \L^m) ) \cr \uparrow & & \uparrow
& & \uparrow \cr
H^2(X,B \otimes \L^m \otimes \sigma_*(-mE)) & \La & H^2(X,B \otimes \L^m) &
\La & H^2(X,B \otimes \L^m \otimes \O_{A_m} ) \cr } $$
Note that
$$
H^2(X,B \otimes \L^m \otimes \O_{A_m} ) = 0 \eqno (2) $$}%
since dim $A_m \leq 1.$

\sn We claim that the pull-back map
$$
\alpha : H^2(X,B \otimes \L^m) \La H^2(\X,\sigma^*(B \otimes \L^m)) $$
is injective. In fact, let $\tau : \tilde X \La \X$ be a desingularisation.
Since $R^1(\sigma \circ \tau)_*(\O_{\tilde X}) = 0,$ the Leray spectral
sequence yields the injectivity of the pull-back map $$\beta : H^2(X,B
\otimes \L^m) \La H^2(\tilde X, (\sigma \circ \tau)^*(B \otimes \L^m)).$$
Since $\beta $ factors through $\alpha,$ also $\alpha$ must be injective.
Now the diagram together with (1) and (2) gives our claim (2.2).

\bn {\bf 2.3 Corollary}
\sn {\it Let $X$ be a compact threefold with $b_2(X) = 0$ and with $a(X) =
2.$ Let $B$ be a vector bundle on $X.$ Then \sn (a) $ H^2(X,B \otimes \L) =
0$ for $\L \in \Pic^0(X)$ generic. \sn (b) $\chi(X,B) \leq 0$ \sn (c)
$c_3(X) = 0.$}

\bn {\bf Proof.} Since $b_2(X) = 0,$ some multiple $\L^{\otimes k}$ with
$\L$ as in (2.2) is in $\Pic^0(X).$ Apply (2.2) to see that $H^2(X,B
\otimes \L^{mk}) = 0$ for large $m$. Hence by semi-continuity (a) holds. By
Serre duality $H^1(X,B \otimes \L) = 0$ for generic $\L.$ From (1.3) we
therefore deduce $\chi(X,B \otimes \L) = 0.$ Hence $\chi(X,B) = 0.$ Now (c)
follows from Riemann-Roch for $B = T_X.$ 

\bn {\bf (2.4)} We now treat the case that $d = 1$ and that $f$ is not
holomorphic. Again we fix an ample divisor $A$ on $V$ and let $\L =
\sigma_*(\f^*(\O_V(A))^{**}.$ We will prove :

\bn {\bf 2.5 Proposition}
\sn {\it Let $X$ be a smooth compact threefold with $a(X) =1$, $b_2(X) = 0$
and with non-holomorphic algebraic reduction $f : X \dasharrow V.$ Then
$$
\chi(X,B) \leq 0
$$
for every vector bundle $B$ on $X$.}

\bn {\bf Proof.} As in (2.3) it is sufficient to prove $$
H^q(X,B \otimes \H) = 0
$$
for $q = 0$ and $q = 2$ and some $\H \in \Pic^0(X)$ (which of course can be
taken individually in each case). \sn For $q = 0$ the claim follows from
(1.3). Turning to $q = 2,$ we fix a positive integer $k$ such that
$$
H^2(\hat F,\sigma^*(B \otimes \L^k)) = 0 \eqno (*), $$
where $\hat F$ is the general fiber of $\hat f.$ This $k$ can be found since $$
H^2(\hat F,\sigma^*(B \otimes \L^k)) =
H^0(\hat F,\sigma^*(B^* \otimes \L^{-k}) \otimes \omega^{-1}_{\hat F}), $$
which vanishes by (1.1), having in mind that $\sigma^*(\L)\vert \hat F =
\hat f^*(A) \otimes\O_{\hat X}(E) \hat F $ is effective and non-zero.

Now let
$$
\H = \L^{m+k}
$$
with $m$ large to be specified in a moment. Arguing as in (2.2) we must
prove $$ H^2(\X,\sigma^*(B \otimes \L^{m+k})(-mE)) = H^2(\X,\sigma^*(B
\otimes \L^k) \otimes \f^*(mA)) = 0. $$
By the Leray spectral sequence for $\f$ we need for that $$
H^1(V,R^1\f_*(\sigma^*(B \otimes \L^k)) \otimes mA) = 0 $$
and
$$
H^0(V,R^2\f_*(\sigma^*(B \otimes \L^k)) \otimes mA) = 0. $$
The first is just Serre vanishing; for the second we claim $$
R^2\f_*(\sigma^*(B\otimes \L^k)) = 0.
$$
But this is just (*) together with (1.5), saying that $R^2\f_*(\sigma^*(B
\otimes \L^k))$ is locally free. This finishes the proof.

\section{3. The case of a holomorphic algebraic reduction}

\bn {\bf 3.0 Setup}

\sn In this section we fix a smooth compact 3-fold $X$ with {\it
holomorphic} algebraic reduction $f : X \La V$ to the smooth curve $V$ of
genus $g.$ We shall always assume $b_2(X) = 0.$ Of course $f$ has connected
fibers.

\bn {\bf 3.1 Theorem}
\sn {\it Let $B$ be a holomorphic vector bundle on $X.$ Then $\chi(X,B) = 0.$}

\bn {\bf 3.2 Corollary }
\sn {\it $c_3(X) = 0.$}

\bn {\bf Proof.} This is just Riemann-Roch for $\chi(X,T_X).$

\bn {\bf Proof of 3.1} \ Let $W = H^1(X,\O_X).$ Then every element in $W$
is represented as a topo\-logically trivial line bundle. Let $A$ be an
ample line bundle on $V.$ Then for $m \gg 0$ we have

$$
H^0(V,R^jf_*(B) \otimes A^{-m}) = 0 \eqno (1)$$ for all $j \geq 0.$ We
claim that (1) implies
$$
H^0(V,R^jf_*(B \otimes \L)) = 0 \eqno (2) $$ for $\L \in W$ general. In
fact, consider locally the universal bundle $\hat {\L}$ on $X \times W.$
Let $F = f \times \id : X \times P \La V \times P$ and $\hat B =
\pr_X^*(B).$ The coherent sheaf $R^jF_*(\hat B \otimes \hat {\L}) $
satisfies
$$
R^jF_*(\hat B \otimes \hat {\L}) \vert V \times \{t\} \simeq R^jf_*(B
\otimes \hat {\L}_t), $$
where $\hat {\L}_t$ is the line bundle corresponding to $t \in W.$ Choose
$m \gg 0$ and $t_0$ such that $f^*(A^{-m}) = \hat \L_{t_0}.$ By (1) and the
assumption $b_2(X) = 0$ we have $$H^0(V,R^jf_*(B \otimes \hat {\L}_{t_0}))
= 0.$$ Hence it is sufficient to show that $R^jF_*(\hat B \otimes \hat
{\L}) $ is flat with respect to the projection $ q :V \times P \La W,$ over
a Zariski open set of $W$, then the usual semi-continuity theorem gives the
claim (2). Now $R^jf_*(B \otimes \hat {\L}_t) $ is locally free on $V = V
\times t$ for every $t$ by (1.5), hence it is clear that there is a Zariski
open set $U \subset W$ such that $R^jF_*(\hat B \otimes \hat {\L})$ has
constant rank over $U$, hence is locally free over $U$ (observe just that
the set where the rank of a coherent sheaf is not minimal is analytic).

\sn On the other hand we have for $m \gg 0 $ and all $i \geq 1, j \geq 0$
by Serre's vanishing theorem
$$
H^{i}(V, R^jf_*(B) \otimes A^m) = 0 \eqno (3).$$ In the same way as (2)
follows from (1) we conclude from (3) that $$
H^{i}(V, R^jf_*(B\otimes {\L})) = 0 \eqno(4) $$
for general $\L \in P.$

\sn Now (2) and (4) give by the Leray spectral sequence $H^{i}(X, B \otimes
\L) = 0$ for general $\L \in P$ and all $i \geq 0.$ In particular $\chi(X,
B \otimes \L) = 0.$ Since $c_k(B \otimes \L) = c_k(B)$, Riemann-Roch gives
our claim $\chi(X, B) = 0.$ 

\section{4. The Main Theorem and Examples}

\bn Putting (2.6), (2.8) and (3.2) together, we obtain

\bn {\bf 4.1 Theorem}
\sn {\it Let $X$ be a compact complex manifold of dimension $3.$ If $b_2(X)
= 0$ and if $c_3(X) \ne 0,$ then every meromorphic function on $X$ is
constant.}

\bn {\bf Proof.} We only have to notice that $a(X) < 3$ by (1.2) and then
to apply the above mentioned propositions.

\bn {\bf 4.2 Corollary}
\sn {\it Let $X$ be a compact complex threefold homeomorphic to the sphere
$S^6.$ Then every meromorphic function on $X$ is constant.}

\bn {\bf Proof.} Note that $c_3(X) = \chi_{\rm top}(S^6) = 2$ and apply (4.1).

\bn The condition $c_3(X) \ne 0$ in (4.1) can be translated into $2b_1(X) +
b_3(X) \ne 2.$

\bn {\bf 4.3 Corollary}
\sn {\it Let $X$ be a compact threefold with $b_2(X) = 0$ and $a(X) = 1 $
or $2.$ Let $X \dasharrow V$ be an algebraic reduction to the smooth
variety $V.$ \sn (1) If $\dim V = 1,$ then $V = \bP_1.$ \sn (2) If $\dim V
= 2$ and $\kappa(X) = - \infty,$ then $V$ is rational.}

\bn {\bf Proof.} By the Main Theorem we have $b_1(X) \leq 1.$ Let $\hat X
\La X$ be a sequence of blow-ups with smooth centers such that the induced
map $ \hat f : \hat X \La V$ is holomorphic. Notice $b_1(\hat X ) \leq 1.$
Now the Leray spectral sequence gives $b_1(V) \leq b_1(\hat X),$ since
$\hat f$ has connected fibers and therefore $\hat f_*(\bC) = \bC.$ So $V =
\bP_1$ in case dim $V = 1.$ \sn If dim $V = 2$ we still need to know that
$\kappa(V) = - \infty .$ This comes from [Ue87].

\bn In the rest of this section we give examples of threefolds with $a(X) > 0$, $b_2(X) = 0$ and $c_3(X) = 0$ so that (4.1) is sharp.

\bn {\bf 4.4 Example}
\sn The so-called Calabi-Eckmann threefolds are compact threefolds
homeomorphic to $S^3 \times S^3,$ see [Ue75]. They can be realised as
elliptic fiber bundles over $\bP_1 \times \bP_1.$ Hence $a(X) = 2$, $b_1(X)
= b_2(X) = 0$ and $b_3(X) = 2.$

\sn We now show that Calabi-Eckmann manifolds can be deformed to achieve
$a(X) = 1 $ or $a(X) = 0$ and $b_1 = b_2 = 0$, $b_3 = 2.$ We choose
positive real numbers $a,b,c$ and let $B = \bC^2 \setminus \{(0,0)\}.$ We
define the following action of $\bC$ on $B \times B :$ $$ (t,x,y,u,v)
\mapsto (\exp(t)x,\exp(at)y,\exp(ibt)u,\exp(ict)v).$$ One checks easily
that this action is holomorphic, free and almost proper so that the
quotient $X$ exists and is a compact manifold. If $a = 1$ and $b = c$, then
$X$ is a Calabi-Eckmann manifold. If however $\log a \ne \bQ$ and $b = c$
resp. $\log a \ne \bQ$ and $\log b, \log c$ are $\bQ$-linearly independent,
then $a(X) = 1$ resp. $a(X) = 0.$

\bn {\bf 4.5 Example}
\sn Hopf threefolds of the form
$$
\bC^3 \setminus \{(0,0,0)\} /\bZ ,
$$
with the action of $\bZ\simeq\{\lambda^k\,;\,k\in\bZ\}$ being defined by $$
\lambda (x,y,z) = (\alpha x, \beta y, \gamma z),\qquad 0 < \vert \alpha
\vert, \vert \beta \vert, \vert \gamma \vert < 1 $$
are homeomorphic to $S^1 \times S^5.$ They have $a(X) = 0, 1$ or $2$ and
$b_1(X) = 1,$ while $b_2(X) = b_3(X) = 0.$ This realises the other
possibility for the pair $(b_1,b_3)$ when $b_2 = 0$ and $a(X) > 0,$ see the
Main Theorem.

\bn Notice that the algebraic reduction is holomorphic in (4.4) but not in
(4.5).

\bn {\bf 4.6 Example }
\sn We finally give other examples of compact threefolds $X$ with $a(X) =
0$ and $b_1 = b_2 = 0$, $b_3 = 2.$ Let $\Gamma \subset \Sl(2,\bC)$ be a
torsion free cocompact lattice in such a way that $X := \Sl(2,\bC)/\Gamma$
has $b_1(X) = 0.$

\sn This last condition is not automatic. Let $Y = \SU(2)\backslash
\Sl(2,\bC)/\Gamma;$ then $Y$ is a compact differentiable manifold admitting
a differentiable fibration $\pi: X \La Y$ with $S^3$ is fiber. Since
$b_1(X) = 0,$ we also have $b_1(Y) = 0,$ hence $b_2(Y) = 0$ by Poincar\'e
duality. Now the Leray spectral sequence immediately gives $b_2(X) = 0,$
the fibers of $\pi$ being 3-spheres. $b_3(X) = 2$ is again clear from the
Leray spectral sequence. Finally the fact that $X$ does not carry any
non-constant meromorphic function results from [HM83] from which we even
deduce that $X$ does not carry any hypersurface (anyway $X$ is
homogeneous). 

\section{5. On the finer structure of threefolds with
$b_2(X)=0$ and $a(X)=1.$}

\bn \bn In this section we investigate more closely threefolds $X$ with
$b_2(X) = 0$ and holomorphic algebraic reduction $f : X \La V$ to a curve
$V.$ We know already that $V$ is rational. We fix these notations. The key
to our investigations is

\bn {\bf 5.1 Theorem}
\sn {\it The restriction $r : H^1(X,\O_X) \La H^1(F,\O_F) $ is surjective
where $F$ is the general smooth fiber of $f.$}

\bn We need some preparations for the proof of (5.1). Let $\Delta \subset
V$ be a finite {\bf non empty } set such that $A = f^{-1}(\Delta) \subset
X$ contains all singular fibers of $f.$ Let $W = V \setminus \Delta, U =
f^{-1}(W)$ so that $g = f \vert U$ has only smooth fibers, hence is a
$C^{\infty}$-fibration. Let $D_i, \ 1\leq i \leq r$ be the irreducible
components of $A$ and let $s = {\rm card} \ \Delta$ be the number of
connected components of $A.$ Furthermore we set $t = b_1(F), \ F$ the
general smooth fiber of $f$ and recall that $g$ is the genus of $V.$ For a
non-compact space $Z$ we let $b_i(Z) = \dim H_i(Z,\bR),$ if this space is
finite-dimensional. We prepare the proof of (3.1) by three lemmas.

\bn {\bf 5.2 Lemma}
\sn {\it (1) The natural exact sequence of groups $$
1 = \pi_2(W) \La \pi_1(F) \La \pi_1(U) \La \pi_1(W) \La 1 $$
is exact and (non-canonically) split.
\sn (2) $b_1(U) = b_1(W) + b_1(F) = 2g + s-1 + t \leq b_1(X) + s-1 +t.$ \sn
(3) $r = s-1 + t + 2g - b_1(X).$}

\bn {\bf Proof.} (1) Since $W$ is a non-compact Riemann surface, $\pi_1(W)$
is a free group of $2g+s-1$ generators and since $W$ is uniformised by
either $\bC$ or by the unit disc, we have $\pi_2(W) = 1.$ Hence the exact
homotopy sequence of the fibration $g : U \La W$ gives the exact sequence
of groups stated in (1). Since $\pi_1(W)$ is a free group, the sequence
splits.

\sn (2) From (1) we deduce that
$$
H_1(U,\bZ) \simeq H_1(W,\bZ) \oplus H_1(F,\bZ).$$ Moreover $f_* : \pi_1(X)
\La \pi_1(V)$ is surjective since the fibers of $f$ are connected. Hence
(2) follows.

\sn (3) The cohomology sequence with rational coefficients of the pair
$(X,A)$ gives
$$
0 = H^4(X) \La H^4(A) \La H^5(X,A) \La H^5(X) \La H^5(A) = 0. $$
By duality we have $H^5(X,A) \simeq H_1(U)$ and $H^5(X) \simeq H_1(X).$
Hence 
$$r = \dim H^4(A) = b_1(U) - b_1(X) = 2g + s-1+t - b_1(X)$$ by (2), as
claimed.

\bn Now choose an integer $m > 0$ such that $c_1(\O_X(mD_i)) = 0$ in
$H^2(X,\bZ)$ for all $1 \leq i \leq r.$ Let $G = \bigoplus \bZ[mD_i] \simeq
\bZ^r$ be the free abelian group of $r$ generators and let $\phi: G \La
\Pic^0(X)$ be given by sending $D = \sum _i a_imD_i$ to $\O_X(D).$

\bn {\bf 5.3 Lemma}
\sn {\it Let $K = {\Ker} \phi$ and $H = {\Im} \phi.$ Then $\rk K \leq s-1 $
and $\rk H \geq r-s+1.$ }

\bn {\bf Proof.} Write $\Delta = \{x_1, \ldots , x_s\}$ and define $G_0 =
\bigoplus _j \bZ[mx_j] \subset {\rm Div}(V).$ Let $K_0 = \{g \in G \vert
{\rm deg}(g) = 0 \}.$ Then $K_0 \simeq \bZ^{s-1}.$ Consider the pull-back
map $f^* : {\rm Div}(V) \La {\rm Div}(X)$ and let $K_1 = f^*(K_0).$ We
claim that $K \subset K_1,$ which will prove the lemma, since then $\rk K
\leq \rk K_1 \leq \rk K_0 \leq s-1.$ To show $K \subset K_1,$ consider $D =
\sum _i a_imD_i = \sum \alpha _j mD_j - \sum \beta_k mD_k = D' - D'' $ with
$\alpha_j, \beta_k > 0$ and $D_j \ne D_k$ for all $j \ne k.$ If $D \in K,$
then $\O_X(D) \simeq \O_X,$ hence $\O_X(D') \simeq \O_X(D'').$ Letting $\L
= \O_X(D'),$ we see that $h^0(X,\L) \geq 2,$ therefore $\L \in f^*(\Pic(V)$
since $f$ is the algebraic reduction of $X.$ Thus $D_i \in K_1,$ hence $D
\in K_1,$ proving the claim.

\bn The last ingredient in the proof of (5.1) is provided by

\bn {\bf 5.4 Lemma}
\sn {\it The kernel of the restriction map $\alpha: H^1(X,\O_X) \La
H^1(U,\O_U)$ equals the image of $f^* : H^1(V,\O_V) \La H^1(X,\O_X).$}

\bn {\bf Proof.} The Leray spectral sequence gives the following
commutative diagram with exact rows :
$$
\matrix{
0 \La & H^1(V,\O_V) \buildrel {f^*} \over \La & H^1(X,\O_X) \La &
H^0(V,R^1f_*(\O_X)) \La 0 \cr
& \downarrow & \downarrow & \downarrow & \ \cr 0 = & H^1(W,\O_W) \La &
H^1(U,\O_U) \La & H^0(W,R^1g_*(\O_W)) \La 0 \cr}
$$
Now $R^1f_*(\O_X)$ is locally free by (1.5). This implies that the
restriction map
$$
\tilde \alpha : H^0(V, R^1f_*(\O_X)) \La H^0(W, R^1g_*(\O_W)) $$
is surjective, hence also $\alpha : H^1(X,\O_X) \La H^1(U,\O_U)$ by the
diagram. This proves (5.4).

\bn We are now able to finish the proof of (5.1).

\sn We consider again $H \subset \Pic^0(X),$ the image of $ \phi : G \La
\Pic^0(X).$ $H$ has rank $\geq r-s+1 = t+(2g-b_1(X))$ by (5.2). Let $\tilde
H$ be the inverse image of $H$ under the natural map $H^1(X,\O_X) \La
\Pic^0(X).$ Then
$$
\rk \tilde H = \rk H + b_1(X) \geq 2g+t. $$
Let $H_1 = \alpha (\tilde H) \subset H^1(U,\O_U).$ We claim that $$
\rk H_1 \geq t. \eqno (*)
$$
In fact, by (5.4) $ \tilde H_0 := {\Ker} (\alpha \vert \tilde H) = \tilde H
\cap f^*(H^1(V,\O_V)).$ Now $\tilde H_0$ is a discrete subgroup of
$f^*(H^1(V,\O_V) \simeq \bC^g,$ so that $rk \tilde H_0 \leq 2g.$ Hence $\rk
H_1 \geq 2g+t - 2g = t,$ proving (*).

\sn Let $\zeta: H^1(U,\bZ) \La H^1(U,\O_U)$ be the canonical map; then $H_1
\subset \zeta(H^1(U,\bZ))$ since $\O_U(D) \simeq \O_U$ for all $D \in G.$
We therefore conclude that the group $H_2 := \beta (H_1)$ with $$\beta :
H^1(U,\O_U) \La H^0(W, R^1g_*(\O_U))$$ the canonical isomorphism, is of
finite index in
$$
\Gamma := {\Im} (H^0(W, R^1f_*\bZ)) \La H^0(W, R^1g_*(\O_U)). $$
By 5.2(2) we conclude that $\rk H^0(W, R^1f_*(\bZ)) = t,$ hence $H_2$ has
finite index in $\Gamma.$ Now take a general smooth fiber $F$ and let
$\lambda : H^0(W, R^1g_*(\O_U)) \La H^1(F,\O_F)$ be the restriction map.
Then we conclude that $\lambda(H_2)$ is of finite index in
$${\Im}(H^1(F,\bZ) \La H^1(F,\O_F)).$$ Since the linear span of this image in
$H^1(F,\O_F)$ is $H^1(F,\O_F),$ we deduce from the construction of $H_2$
that the restriction map $H^1(X,\O_X) \La H^1(F,\O_F)$ is surjective
proving finally (5.1).

\bn {\bf (5.5)} We now study the structure of the smooth fibers $F$ of $f.$
Since $K_F \equiv 0$ and $c_2(F) = c_2(X) \vert F = 0$, we conclude from
the classification of surfaces that $F$ is one of the following: a Hopf
surface, an Inoue surface, a Kodaira surface (primary or secondary), a
torus or hyperelliptic (see e.g. [BPV84]). By [Ka69] however, $F$ cannot be
hyperelliptic. The reason is the existence of a relative Albanese reduction
in that case.

\bn {\bf 5.6 Theorem}
\sn {\it The general fiber $F$ of $f$ cannot be a Kodaira surface nor a
Hopf surface with algebraic dimension $1.$ }

\bn {\bf Proof.} Let $F_0$ be a fixed smooth fiber and assume that $F_0$ is
a Kodaira surface or a Hopf surface with $a(F_0) = 1.$ Let $g_0 : F_0 \La
C_0$ be "the" algebraic reduction which is an elliptic fiber bundle. Let
$\L_0 = g_0^*({\cal G}_0)$ with ${\cal G}_0$ very ample on $C_0.$

\sn (1) There exists a line bundle $\tilde \L$ on $X$ with $\tilde \L \vert
F_0 = \L_0.$ \sn Proof. By passing to some power $\L_0^m$ if necessary, we
have $c_1(\L_0) = 0$ in $H^2(X,\bZ).$ Let $$
\lambda_1 : H^1(X,\O_X) \La {\rm Pic}(X) $$
and
$$
\lambda_2 : H^1(F_0, \O_{F_0}) \La {\rm Pic}(F_0) $$
be the canonical maps, $r : H^1(X,\O_X) \La H^1(F_0,\O_{F_0}) $ the
restriction map. Choose $\alpha \in H^1(F_0,\O_{F_0})$ with
$\lambda_2(\alpha) = \L_0.$ Since $r$ is surjective by (4.2), we find
$\beta \in H^1(X,\O_X)$ with $r(\beta) = \alpha.$ Now let $\tilde \L =
\lambda_1(\beta).$ \sn (2) Let $F$ be any smooth fiber of $f.$ Then
$\kappa(\tilde \L \vert F) = 1.$ In fact, it follows from the local
freeness of $R^jf_*(\tilde \L)$ (1.5) that $$
f_*(\tilde \L^{\mu}) \vert \{y\} \simeq H^0(F,\tilde \L^{\mu}), $$
where $F = f^{-1}(y),$ see [BaSt76, chap.3, 3.10]. \sn (3) From the
generically surjective morphism
$$
f^*f_*(\tilde \L^m) \La \tilde \L^m,
$$
($m \gg 0)$, we obtain a meromorphic map $$
g: X \dasharrow \bP(f_*(\tilde \L^m)),$$ which, restricted to $F$ is
holomorphic and just gives the algebraic reduction of $F.$ Let $Z$ be the
closure of the image of $g.$ Then $f$ factors via the meromorphic map $ h_1
: X \dasharrow Z$ and the holomorphic map $h_2 : Z \La V.$ Now $h_2$ is the
restriction of the canonical projection $\bP(f_*(\tilde \L^m)) \La V,$
therefore $h_2$ is a projective morphism and $Z$ is projective. Hence
\hbox{$a(X) \geq 2,$} contradiction.

\bn From (5.6) it follows that $F$ can only be an Inoue surface, a Hopf
surface without meromorphic functions or a torus. In order to exclude by a
similar method as in (5.6) also tori of algebraic dimension 1, we would
need the existence of a relative algebraic reduction (the analogue of $h_2
: X \dasharrow Z$) in that case, too.

\sn We now look more closely to the structure of $f.$

\goodbreak
\bn {\bf 5.7 Proposition}
\par\nobreak
\sn {\it Assume that $F$ is not a torus. Then \item {(1)} $R^1f_*(\O_X) = \O_V $
\item {(2)} $R^2f_*(\O_X) = 0 $
\item {(3)} ${\dim} H^1(X,\O_X) = 1$
\item {(4)} $H^2(X,\O_X) = H^3(X,\O_X) = 0.$ }

\bn {\bf Proof.} (2) Since $H^2(F,\O_F) = 0,$ the sheaf $R^2f_*(\O_X) $ is
torsion, hence $0$ by (1.5). Then $H^3(X,\O_X) = 0$ is immediate from the
Leray spectral sequence.

\sn (1) Since $h^1(\O_F) = 1, R^1f_*(\O_X)$ is a line bundle on $V.$ Let
$$
d = {\rm deg} R^1f_*(\O_X).
$$
Then Riemann-Roch gives $\chi(R^1f_*(\O_X)) = d+1. $ On the other hand the
Leray spectral sequence together with $H^3(X,\O_X) = 0$ and (2) yields
$$
\chi(R^1f_*(\O_X)) = h^1(\O_X) - h^2(\O_X) = - \chi(\O_X) + 1. $$
We conclude $d = - \chi(\O_X) = 0.$ This proves (1). Now (3) and the second
part of (4) are obvious.

\bn {\bf 5.8 Remark} \ In case $F$ is a torus, $R^1f_*(\O_X)$ is a rank 2
bundleand $R^2f_*(\O_X)$ is a line bundle. Using (4.2) it is easy to see
that {\item {(1)} $R^1f_*(\O_X) = \O(a) \oplus \O(b) $ with $a,b \geq 0$}
{\item {(2)} $R^2f_*(\O_X) = \O(a+b).$}

\sn Note that (2) gives dually $f_*(\omega_{X\vert V}) = \O(-a-b)$. Usually
one expects the degree of $f_*(\omega_{X \vert V})$ to be semi-positive,
but here we are in a highly non-K\"ahler situation where it might happen
that the above degree is negative, see [Ue87].

\bn {\bf 5.9 Proposition}
\sn {\it Assume that $F$ is not a torus. Then $H^0(X,\Omega^{i}_X) = 0$, $1
\leq i \leq 3.$}

\bn {\bf Proof.} For $i = 3$ the claim follows already from (5.7) and (5.8).

\sn (1) First we treat the case $i = 1.$ Let $\omega$ be a holomorphic
1-form. Let $j : F \La X$ be the inclusion. Then $j^*(\omega) = 0$, hence
at least locally near $F$ we have $\omega = f^*(\eta),$ hence $d\omega = 0$
near $F$ and therefore the holomorphic 2-form $d\omega = \partial \omega =
0$ everywhere. Now the space of closed holomorphic 1-forms can be
identified with $H^0(X,d\O_X)$ and as it is well known (see e.g. [Ue75]),
we have the inequality $$ 2 h^0(X,d\O_X) \leq b_1(X),
$$
hence $b_1(X) \leq 1 $ gives our claim.

\sn (2) In case $i = 2,$ we again have $j^*(\omega) = 0.$ Let $U$ be a
small open set in $V$ such that $f \vert f^{-1}(U)$ is smooth. $z$ be a
coordinate on $U$ and $h = f^*(z).$ Then we conclude that $$
\omega \vert f^{-1}(U) = dh \wedge\alpha $$
with some holomorphic 1-form $\alpha.$ Now again $j^*(\alpha) = 0$ and
therefore $\alpha = f^*(\beta).$ In total $\omega = f^*(dz \wedge \beta) =
0.$

\bn {\bf 5.10 Corollary }
\sn {\it Assume that $F$ is not a torus. Then there exists some $x$ such
that $f_*(\Omega^1_X) = \bC_x,$ i.e. a sheaf supported on $x$ with a
1-dimensional stalk at $x.$ In particular $f$ has at most one normal
singular fiber and such a fiber has exactly one singularity which has
embedding dimension 3.}

\bn {\bf Proof.} The exact sequence
$$
0 \La f^*(\Omega_V) \La \Omega^1_X \La \Omega^1_{X \vert V} \La 0 $$
together with (5.7) and the projection formula leads to an exact sequence
$$ 0 \La \Omega^1_V \La f_*(\Omega^1_X) \La f_*(\Omega^1_{X \vert V}
{\buildrel \alpha \over \longrightarrow} \ \Omega^1_V \otimes R^1f_*(\O_X)
\simeq \Omega^1_V.
$$
Since $F$ has no holomorphic 1-forms, $f_*(\Omega^1_{X \vert V})$ is a
torsion sheaf, therefore $\alpha = 0.$ Now observe the following facts.
\item {(1) } $ H^0(X,\Omega^1_X ) = 0$ by (5.9) \item {(2) }
$H^1(X,f^*(\Omega_X)) = H^1(V,\Omega^1_V) $ \item {(3) } $H^1(\Omega^1_V)
\La H^1(X,\Omega^1_X) $ is $0$ (since $b_2(X) = 0).$ \sn Then taking
cohomology of the first sequence, it follows that $$ h^0(X,\Omega^1_{X
\vert V}) = 0,
$$
hence the second sequence gives our claim.

\bn We can say something more about the structure of the singular fibers of $f.$

\bn {\bf 5.11 Proposition}
\sn {\it Assume that $F$ is not a torus. Let $A$ be a union of fibers
containing all singular fibers of $f.$ Let $s = {\rm card}(f(A))$ and $r$
the number of irreducible components of $A.$ Then $r = s,$ i.e. all fibers
of $f$ are irreducible and $b_1(X) = 0.$ }

\bn {\bf Proof}. Since $F$ is an Inoue surface or a Hopf surface, we have
$b_1(F) = 0,$ moreover $V = \bP_1,$ hence the claim follows from (5.2).

\bn {\bf 5.12 Remark} In case $F$ is a torus, the same proof gives $r =
s+3.$ It seems rather reasonable to expect that tori actually cannot appear
as fibers of $f.$ Observe that $f$ must have a singular fiber in this case
because of $r = s+3.$ So a study of the singular fibers is needed to
exclude tori as fibers of $f.$ We hope to come back to this question in a
later paper.

\bn {\bf 5.13 Proposition}
\sn {\it Assume that $F$ is not a torus. Then $h^{1,1} = h^{1,2} = h^{2,1}
= 1 $ (so that we know all Hodge numbers of $X.$)}

\bn {\bf Proof.}(1) $h^{1,2} = h^{2,1}$ is of course Serre duality.

\sn (2) By (5.9) and $\chi(X,\Omega_X^1) = 0$ it suffices to see $h^{1,3} =
0$ in order to get $h^{1,1} = h^{1,2}.$ But this follows again from Serre
duality and (5.9).

\sn (3) From the exact sequence
$$
0 \La f^*(\Omega^1_V) \La \Omega^1_X \La \Omega^1_{X \vert V} \La 0 \eqno (S) $$
we deduce that it suffices to show \item {(a)} $h^2(X,f^*(\Omega^1_V)) = 1$
\item {(b)} $h^2(X,\Omega^1_{X \vert V}) = 0$ \sn in order to get $h^{1,2}
\leq 1.$

\sn (a) By the Leray spectral sequence and (5.7) we have
$h^2(X,f^*(\Omega_V^1)) = h^2(V,\Omega^1_V) = 1.$

\sn (b) Again we argue by the Leray spectral sequence. Since
$R^1f_*(\Omega^1_{X \vert V})$ is a torsion sheaf, we need only to show
that
$$
R^2f_*(\Omega^1_{X \vert V}) = 0.
$$
But, taking $f_*$ of (S), $R^2f_*(\Omega^1_{X \vert V})$ is a quotient of
$R^2f_*(\Omega^1_X)$ which is $0$ by $H^2(\O_F) = H^2(\Omega^1_F) = 0$ and
(1.5).

\sn (4) We finally show $h^{1,1} \ne 0$ to conclude the proof. Let
$(E^{pq}_r)$ be the Fr\"olicher spectral sequence on $X.$ Since $b_1(X) =
0$ by (5.11), we get $E_{\infty}^{0,1} = 0.$ Hence $E_2^{0,1} = 0.$ On the
other hand
$$
E_2^{0,1} = {\rm Ker}\partial : E_1^{0,1} \La E_1^{1,1}. $$
Since $E_1^{p,q} = H^{pq}(X),$ we conclude that $H^1(X,\O_X)$ injects into
$H^{1,1}(X).$ So by (5.7) $H^{1,1}(X) \ne 0.$

\bn We finally collect all our knowledge in the case the general fiber of
$f$ is not a torus.

\bn {\bf 5.14 Theorem}
\par\nobreak
\sn {\it Let $X$ be a smooth compact threefold with $b_2(X) = 0$ and
holomorphic algebraic reduction $f: X \La V$ to the smooth curve $V.$
Assume that the general smooth fiber is not a torus. Then: \item
{(1)}$b_1(X) = 0, b_3(X) = 2$
\item {(2)} any smooth fiber of $f$ is a Hopf surface without meromorphic
functions or an Inoue surface. \item {(3)} the Hodge numbers of $X$ are as
follows: $h^{1,0}= 0$, $h^{0,1} = 0$, $h^{2,0} = 0$, $h^{1,1} = 1$,
$h^{0,2} = 0$, $h^{3,0} = 0$, $h^{2,1} = 1$ (the others are determined by
these via Serre duality).
\item {(4)} all fibers of $f$ are irreducible. There is at most one normal
singular fiber. \vskip0pt} 

\ninepoint
\section{References}
\bn

\itemitem{[BaSt76]} Banica, C.; Stanasila, O.: Algebraic methods in the
global theory of complex spaces. Wiley 1976
\vskip2pt
\itemitem{[BPV84]} Barth, W.; Peters, C.; van de Ven,: Compact complex
surfaces. Erg.\ d.\ Math.\ (3), Band 4, Springer 1984
\vskip2pt
\itemitem{[HM83]} Huckleberry, A.T.; Margulis, G.: Invariant analytic
hypersurfaces. Inv.\ Math.\ {\bf 71}, 235-240 (1983)
\vskip2pt
\itemitem{[Ka69]} Kawai, S.: On compact complex analytic manifolds of
complex dimension 3, II. J.\ Math.\ Soc.\ Japan {\bf 21}, 604-616 (1969)
\vskip2pt
\itemitem{[RRV71]} Ramis, J.; Ruget, G.; Verdier, J.L.: Dualit\'e relative en
g\'eom\'etrie analytique complexe. Inv.\ Math.\ {\bf 13}, 261-283 (1971)
\vskip2pt
\itemitem{[Ue75]} Ueno, K.: Classification theory of algebraic varieties and
compact complex spaces. Lecture Notes in Math.\ 439, Springer 1975
\vskip2pt
\itemitem{[Ue87]} Ueno, K.: On compact analytic threefolds with non-trivial
Albanese tori. Math.\ Ann.\ {\bf 278}, 41-70 (1987)
\vskip2pt
\itemitem{[We85]} Wehler, J.: Der relative Dualit\"atssatz f\"ur
Cohen-Macaulay-R\"aume. Schriftenr. des Math.\ Instituts der Univ.\
M\"unster, 2.Serie, Heft {\bf 35} (1985)

\bigskip\rm
\noindent{\ninecsc
Fr\'ed\'eric Campana, D\'epartement de math\'ematiques,
Universit\'e de Nancy, BP 239, 54506 Vandoeuvre les Nancy, France
({\tt frederic.campana@iecn.u-nancy.fr})}

\medskip
\noindent{\ninecsc
Jean-Pierre Demailly, Universit\'e de Grenoble, Institut Fourier, BP 74,
U.R.A.\ 188 du C.N.R.S., 38402 Saint-Martin d'Heres, France
({\tt demailly@ujf-grenoble.fr})}

\medskip
\noindent{\ninecsc
Thomas Peternell, Mathematisches Institut, Universit\"at Bayreuth, 95440
Bayreuth, Germany
({\tt peternel@btm8x1.mat.uni-bayreuth.de})}

\bye